\documentclass[a4paper,10pt]{article}

\usepackage{amsmath,amsfonts,amscd,amssymb,amsthm}
\usepackage[dvips]{graphics}
\usepackage[all]{xy}

\newtheorem{thm}{Theorem}[section]
\newtheorem*{thmm}{Theorem}
\newtheorem*{ThmA}{Convergent statement}
\newtheorem*{ThmB}{Formal statement}
\newtheorem{lemme}[thm]{Lemma}

\newtheorem{definition}{Definition}[section]

\newtheorem{prop}{Proposition}[section]
\newtheorem{propr}{Property}[section]
\newcommand{\dd}{d\hspace{0.04cm}}
\newcommand{\OO}{\mathcal{O}}
\newcommand{\DD}{\mathcal{D}}

\newcommand{\FF}{\mathcal{F}}

\newcommand{\ZZ}{\mathcal{Z}}
\newcommand{\MM}{\mathcal{M}}
\newcommand{\UU}{\mathcal{U}}

\newcommand{\GG}{\mathcal{G}}
\newcommand{\NN}{\mathcal{N}}
\newcommand{\WW}{\mathcal{W}}

\newcommand{\fami}[1]{\left\{ #1 \right\}}
\newcommand{\coup}[1]{\left( #1 \right)}

\newcommand{\carre}{
\begin{flushright}
$\square$
\end{flushright}
\medskip}

\newcommand{\petit}[1]{{\scriptscriptstyle #1}}
\newcommand{\dbar}{\overline{\mathbb{D}}}

\newenvironment{demo}{{\it Proof:}}{\carre}
\newcommand{\fleche}[1]{\stackrel{\scriptscriptstyle #1}{ \xymatrix{\ar@{-->}[r]&  \\ }}}

\begin{document} 

\title{Analytical and formal classifications of quasi-homogeneous foliations in $(\mathbb{C}^2,0)$}
\maketitle

\begin{abstract} 
We prove a result of classification for germs of formal and convergent quasi-homogeneous foliations in $\mathbb{C}^2$ with fixed separatrix. Basically, we prove that the analytical and formal class of such a foliation depend respectively only on the analytical and formal class of its representation of projective holonomy.
\end{abstract}

\section*{Introduction and main statements.}
A \emph{germ of foliation $\FF$ with isolated singularity in $\mathbb{C}^2$} is given by an holomorphic $1$-form up to unity
$$ \omega=a(x,y)\dd x +b(x,y)\dd y $$ where $a,b \in \mathbb{C}\{x,y\}$, $a(0,0)=b(0,0)=0$ and $\textup{gcd}(a,b)=1$. The foliation $\FF$ is said to be \emph{formal}, and then denoted by $\widehat{\FF}$ when $a$ and $b$ are only formal functions in $\mathbb{C}[[x,y]]$. A \emph{separatrix} of $\FF$ is a leaf of the regular foliation given by $\omega|_{({\mathbb{C}^2}^*,0)}$ whose closure in $(\mathbb{C}^2,0)$ is an irreducible analytical subset. A \emph{formal separatrix} of $\widehat{\FF}$ is a formal irreducible curve $\{\widehat{f}=0\}, \widehat{f} \in\mathbb{C}[[x,y]]$ such that $\widehat{f}$ divides the product $\dd \widehat{f}\wedge \omega$. Obviously, any convergent separatrix seen as a formal curve is a formal separatrix of $\FF$ seen as a formal foliation. C. Camacho and P. Sad show in \cite{separatrice} that $\FF$ has at least one separatrix. Once the foliation is desingularized by a blowing-up morphism \cite{reduction} $$E:(\MM,\DD)\rightarrow (\mathbb{C}^2,0),$$ the foliation $\FF$ is said to be of \emph{generalized curve type} when the pull-back foliation $E^*\FF$ has no singularity of \emph{saddle-node} type. C. Camacho, A. Lins Neto and P. Sad prove in \cite{camacho} that the foliation $\FF$ and its separatrix have the same desingularization under the generalized curve assumption. All these definitions and results have their equivalent in the formal context: one just adds the word \emph{formal} in order to get adapted statements.

\noindent According to \cite{MatQuasi}, the foliation $\FF$ is said to be \emph{quasi-homogeneous} when the union  of its separatrix, denoted by $\textup{Sep}(\FF)$, is a germ of curve given in some coordinates by a quasi-homogeneous polynomial function:
$$ \textup{Sep}(\FF)=\left\{\sum_{\alpha i + \beta j = \gamma} a_{ij}x^iy^j=0\right\},\quad \alpha,\beta,\gamma \in \mathbb{N}^*,\quad \textup{gcd}(\alpha, \beta) =1.$$
The couple $(\alpha,\beta)$ is called the \emph{weight} of the curve.. In \cite{MatQuasi}, one can find the following criterion:
\begin{thmm}[J.F. Mattei]
Let $\omega=a(x,y)\dd x + b(x,y) \dd y$ be a germ of holomorphic $1$-form in $\mathbb{C}^2$ with an isolated singularity at $0$. Let us suppose that $\omega$ defines a foliation $\FF$ of generalized curve type with a finite number of separatrix ( i.e. non-dicritical ) and denote by $f(x,y)$ an equation of the separatrix. The following properties are equivalent:
\begin{enumerate}
\item $\FF$ is a quasi-homogeneous foliation.
\item $f$ belongs to the ideal $(a,b)\subset \OO_{\mathbb{C}^2,0}$.
\item $f$ belongs to its jacobian ideal, $(f'_x,f'_y)\subset \OO_{\mathbb{C}^2,0}$.
\item there exists coordinates $u,v$ such that $f$ is written as a\hyphenation{qua-si-ho-mo-ge-neous} quasi-homogeneous polynomial function with $ (\alpha,\beta)$ weight and there exists holomorphic functions $g,h$ with $g(0)\neq 0$ such that
$$g\omega=\dd f +h(\beta u\dd v-\alpha v\dd u).$$
The above expression is called the \emph{normal form of Takens} of the foliation.
\item $\FF$ satisfies the following equivalence: any integrable isoholonomic deformation of $\FF$ is analytically trivial if and only if the underlying deformation of separatrix is analytically trivial.
\end{enumerate}
We refer to the latter property by saying that $\FF$ is \emph{D-quasi-homogeneous}.
\end{thmm}
\smallskip
\noindent From now on, let us suppose $\FF$ to be a quasi-homogeneous foliation of\hyphenation{ge-ne-ra-lized} generalized curve type. The exceptional divisor $\DD$ comes to be a unique chain of irreducible components such that each component meets exactly two others except for the two extremal components. We will see that there is a unique component of the divisor, which meets the closure of $E^{-1}(\textup{Sep}(\FF)\backslash \{0\})$ in $\MM$. In this article, we call the latter component the \emph{central component}. As we will explain, the whole transversal structure of $\FF$ is concentrated in the projective holonomy representation over the central component. 

\smallskip
\noindent Let $\FF_0$ and $\FF_1$ be two germs of quasi-homogeneous foliations of generalized curve type with analytically conjugated separatrix. Let us denote by $\Phi$ a conjugacy of the separatrix and $D_0$ and $D_1$ the respective central components. We say that $\FF_0$ and $\FF_1$ have \emph{same holonomy} if there is an interior automorphism $\phi$ of the group $\textup{Diff}(\mathbb{C},0)$ such that the projective holonomy representations satify a commutative diagram:
$$\begin{CD}
  \Pi_1\left(D_0\backslash \textup{Sing}(E_0^*\FF_0),x\right) @>\textup{Hol}_{\FF_0}>> \textup{Diff}(\mathbb{C},0) \\
  @V\Phi^* VV @V\phi VV\\
    \Pi_1\left(D_1\backslash\textup{Sing}(E_1^*\FF_1),\Phi(x)\right) @>\textup{Hol}_{\FF_1}>> \textup{Diff}(\mathbb{C},0)
\end{CD}
$$
The aim of this article is to show the following theorem:
\begin{ThmA}
Let the foliations $\FF_0$ and $\FF_1$ be two germs of quasi-homogeneous foliations with analytically conjugated separatrix. The foliations $\FF_0$ and $\FF_1$ are analytically conjugated if and only if they have the same\hyphenation{ho-lo-no-my} holonomy.
\end{ThmA}
\noindent In \cite{cermou}, D. Cerveau and R. Moussu generalize a construction of the latter \cite{Moussutoutseul} and prove the above result in the case of a quasi-homogeneous foliation given by a $1$-form with nilpotent linear part
$ y\dd y + \cdots.$ 
Their proof consists in extending the conjugacy of the holonomy in a neighborhood of the divisor deprived of a few points by lifting the path with respect to a transversal fibration as in \cite{MM}. The normal form of Takens allows them to build such a fibration and to establish that their extension is bounded and therefore well defined around the whole divisor. In order to show the generalized theorem, we use a very different method based upon two elements:  the existence of an integrable isoholonomic deformation with prescribed underlying family of separatrix and the property of D-quasi-homogeneity. 

\medskip
\noindent The normal forms of Takens are convergent but, unfortunately, fail to be unique. In \cite{zola}, H. Zoladeck and E. Strozyna build up a unique formal normal form for any formal quasi-homogeneous vector fields in the nilpotent case. In \cite{prepafrank}, F. Loray gives a geometric approach in order to obtain these normal forms. More recently, E. Paul generalized these results and found unique formal normal forms for any formal quasi-homogeneous vector fields \cite{paulnormal}. With the work of J.P. Ramis, M. Canalis-Durand, R. Sch\"a fke and Y. Sibuya \cite{Cana}, one can even compute the Gevrey degree of these normal forms. 
The notion of projective holonomy representation has a natural extension in the formal context. Hence, we will show the following result:
\begin{ThmB}
Let the foliations $\widehat{\FF}_0$ and $\widehat{\FF}_1$ be two germs of formal quasi-homogeneous foliations with formally conjugated separatrix. The foliations $\widehat{\FF}_0$ and $\widehat{\FF}_1$ are formally conjugated if and if only they have the same formal holonomy.
\end{ThmB}
\noindent This theorem ensures the equivalence between the data of a formal normal form and the choice of a point in the space of representation of a free group in $\widehat{\text{Diff}}(\mathbb{C},0)$. However, the meaning of such a correspondance is still to be worked out. 

\section{Isoholonomic deformations with prescribed\hyphenation{se-pa-ra-trix} separatrix.}
When one wants to use the method of R. Moussu and D. Cerveau in the general case, one has some trouble to extend the conjugacy of foliation. Basically, the obstruction comes from the existence of non quasi-homogeneous functions in the equisingularity class of a quasi-homogeneous function. In this section, we are going to establish a result of construction of isoholonomic deformations which removes the mentioned obstruction. This result has its own interest since it is a very general one. The reader which would be only interested in the classification result could admit the theorem (\ref{exis}) and begin its lecture with the part \ref{part2}.

\subsection{Isoholonomic deformations.}
Let $\FF$ be a formal foliation in $\mathbb{C}^2$. Let $K$ be a compact connected subset of $\mathbb{C}^p$. In this article, we call \emph{transversally formal integrable deformation over $K$of $\FF$} any formal foliation $\FF_K$ of codimension one in $(\mathbb{C}^{2+p},0\times K)$ with $0\times K$ as singular locus such that the leaves are transversal to the fiber of the projection 
$$\pi:(\mathbb{C}^{2+p},0\times K)\rightarrow (\mathbb{C}^p,K),\quad \pi(x,t)=t.$$ Such a deformation is given by an holomorphic integrable $1$-form 
\begin{eqnarray*}
\Omega(x,y,t)=a(x,y,t)\dd x + b(x,y,t)\dd y +\sum_{i=1}^p c_i(x,y,t)\dd t_i,\\
 (x,y)\in \mathbb{C}^2, t=(t_1,\ldots,t_p)\in K
\end{eqnarray*}
with: $a,b,c_1,\ldots,c_p$ lies in the ring $\OO_K[[x,y]]$ of formal series with coefficients in the ring of holomorphic functions on $K$ ; $\Omega$ is integrable; the set $$\fami{a=0,b=0,c_1=0,\ldots,c_p=0}$$ is equal to $0\times K$ and the ideal $(c_1,\ldots, c_p)$ is a sub-ideal of $\sqrt{(a,b)}$. The latter condition is equivalent to the transversality of leaves and fibers of $\pi$. If the coefficients come to be convergent series the deformation is simply called \emph{integrable deformation}. In any case, the foliation $i_t^{*}\FF$ where $i_t,\  t\in K$ is the embedding $i_t(x)=(x,t),\ x\in \mathbb{C}^2$, defines a formal foliation in usual meaning. More generally, for any subset $J$ of $K$, we denote by $\FF_J$ the tranversally formal integrable deformation over $J$ induced by restriction of $\FF_K$ on $\pi^{-1}(J)$. A transversally formal integrable deformation is said \emph{equisingular} if the induced family of foliation $\fami{i_t^{*}\FF}_{t\in K}$ admits a reduction of singularities \emph{in family}. We refer to \cite{MatSal} for a precise definition. 

\medskip
\noindent In order to make more specific the arguments we will use in the section \ref{findelapreuve}, we recall a classical result about a particular case of germ of equisingular integrable deformation with parameters in $(\mathbb{C},0)$. Let $\FF$ be any covergent germ of foliation in $\mathbb{C}^2$. Let $E:(\MM,\DD)\rightarrow (\mathbb{C}^2,0)$ be the desingularization of $\FF$. Let us denote by $\textup{Fix}(\FF)$ the sheaf over $\DD$ whose fiber is the group of local automorphisms $\phi$ in $(\MM,\DD)\times(\mathbb{C},0)$ such that
\begin{itemize}
\item $\phi|_{\MM\times \{0\}}=\textup{Id}$
\item $\phi$ commutes with the projection $\mathtt{P}:(\MM,\DD)\times(\mathbb{C},0)\rightarrow (\mathbb{C},0)$
\item $\phi$ lets invariant each local leaf and $\phi^*(\FF\times(\mathbb{C},0))=\FF\times(\mathbb{C},0)$ 
\end{itemize}
Basically, the flows of vector field $X$ tangent to the foliation $\FF\times(\mathbb{C},0)$ with $D\mathtt{P}(X)=0$ are sections of $\textup{Fix}(\FF)$.
\begin{thm}[J.F. Mattei]\label{intco}
There is a bijection between the moduli space of germs of equisingular integrable deformations of $\FF$ with parameters in $\mathbb{C}$ and $H^1(\DD,\textup{Fix}(\FF))$.
\end{thm}

\noindent Since the sections of $\textup{Fix}(\FF)$ act in the local leaf, the holonomy pseudo-group of $i_t^*\FF_{(\mathbb{C},0)}$ does not depend on $t$ along any equisingular integrable deformation. Basically, it means that \emph{the complex structure on the spaces of leaves} does not vary. Hence, we adopt the following defintion:
\begin{definition}
An (resp. transversally formal) isoholonomic deformation of foliation is an equisingular (resp. transversally formal) integrable deformation.
\end{definition}

\begin{quote} {\bf For simplicity, the article is written from now on in the convergent context. However, there is no special difficulty for transposing the proofs and the results in the formal context.}
\end{quote}

\subsection{Existence of isoholonomic deformation with prescribed separatrix.}
In order to give a precise statement, we have to introduce some more definitions. When none singularity of the desingularized foliation associated to $i_t^{*}\FF_K$ is \emph{a saddle-node} \cite{MatQuasi}, the deformation $\FF_K$ is said to be of generalized curve type. This definition is coherent since, along an isoholonomic deformation, the property holds for any foliation $i_t^{*}\FF_K$ as soon as it holds for one $t\in K$. The generalized curve type hypothesis is obviously a generical one. A \emph{separatrix} of $\FF_K$ is an invariant hypersurface of the regular foliation $\FF|_{{\mathbb{C}^2 }^*\times K}$ such that its closure is an irreducible analytical germ of hypersurface along $0\times K$. When $\FF_K$ has only a finite number of separatrix, we denote by $\textup{Sep}(\FF_K)$ their union and $\FF_K$ is said \emph{non-dicritical}.

\medskip
\noindent From now on, we assume the compact of parameters $K$ to have a fundamental system of open connected neighborhoods, which are Stein open sets. Let $C$ be an analytical subset of $K$ given by the zeros of holomorphic functions.
We do not assume $C$ to be connected in $K$.
\begin{thm}\label{exis}
Let $\FF^0$ be a non-dicritical isoholonomic deformation over $K$ of generalized curve type and let us denote by $S^0$ the equisingular family $\textup{Sep}(\FF^0)$. Let $S^1$ be any equisingular family over $K$ of germs of curves at the origin of $\mathbb{C}^2$ with
\begin{enumerate}
\item $S^1$ and $S^0$ are topologically equivalent,
\item $S^1|_C$ and $S^0|_C$ are analytically equivalent.
\end{enumerate}
There exists an isoholonomic deformation $\FF^1$ over $K$ such that 
\begin{enumerate}
\item $\textup{Sep}(\FF^1)=S^1,$
\item  $\FF^1|_C$ and $\FF^0|_C$ are analytically equivalent.
\end{enumerate}
Moreover, the deformation $\FF_1$ and $\FF_0$ are embedded in an isoholonomic deformation $\FF$ over $K\times \dbar$ 
$$ \FF|_{K\times\{0\}}=\FF_0\quad \FF|_{K\times\{1\}}=\FF_1$$
which is trivial along any subset $\{t\}\times \dbar,\ t\in C.$
\end{thm}
\noindent Our proof of this theorem is based upon the cohomological interpretation of isoholonomic deformation. Essential arguments can be found in \cite{moi} in the special case $K=\{0\}$. In this section, we are going to follow the proof performed in \cite{moi} and give only the arguments related to the difficulties appearing with the parameter context. One of the interests of the formalism introduced in \cite{moi} is the easy way arguments are generalized to parameter context. The main difficulties are removed thanks to the Stein property of $K$.

\noindent  In the formal context, we establish a result similar to (\ref{exis}). One just has to replace convergent objects by transversally formal objects and the isoholonomic deformation by a transversally formal isoholonomic deformation.

\noindent In the section (\ref{catego}), we fix a manifold $\MM$, which is built over the reduction tree of a isoholonomic deformation and we define a very special class of manifold denoted by $\text{Glu}^C_0(\MM,\UU,Z)$ related to $\MM$. In the section (\ref{demons}), $\MM$ is supposed to be foliated by an isoholonomic deformation $\FF$. A property of cobordism type is pointed out and allows us to detect the existence of an isoholonomic deformation on any element of $\text{Glu}^C_0(\MM,\UU,Z)$. This deformation will automatically be linked to $\FF$ by isoholonomic deformations. In sections (\ref{demo1}), (\ref{demo2}) and (\ref{demo3}), we show the theorem (\ref{cobordisme0}) which states that, under the generical hypothesis, the cobordism property holds for any element of $\text{Glu}^C_0(\MM,\UU,Z)$. In the third section, we deduce the theorem (\ref{exis}) from this cobordism property.

\subsection{The category $\text{Glu}^C_n(\MM,Z,\UU)$.}\label{catego}
A blowing-up process over $K$ is a commutative diagram 
\begin{equation}\label{ProEcla}
\begin{array}{cccccccccl}
\MM^h & \stackrel{E^h}{\rightarrow} & \ldots &\MM^{j} & \stackrel{E^j}{\rightarrow} &  & \ldots & \stackrel{E^1}{\rightarrow} & \MM^{0} & = 
(\mathbb{C}^{2+p},0\times K)\stackrel{\pi}{\rightarrow} (\mathbb{C}^p,K) \\
\bigcup &  &    &\bigcup &  &
 &  &  & \bigcup &  \\
\Sigma^h & \rightarrow & \ldots &\Sigma^{j} & \rightarrow &  & \ldots & \rightarrow & \Sigma^{0} & =  0\times K\\
\bigcup &  &  &  \bigcup &  &  & & & \bigcup &  \\
S^h & \rightarrow & \ldots  &S^{j} & \rightarrow &  & \ldots & \rightarrow & S^{0} & =  0\times K\\
\end{array}\end{equation}
where $\MM^j$ is an analytical manifold of dimension 2+p; $\Sigma^j$ is a closed analytical subset of dimension $p$ called the $j^{th}$  singular locus; $S^j\subset \Sigma^j$ is an analytical smooth sub-variety with a finite number of connected components called the $j^{th}$ blowing-up center. Futhermore, denoting $$E_j:=E^0\circ\cdots \circ E^j,\quad \pi_j:=\pi\circ E_j,\quad \text{et} \quad \DD^j:=E_{j}^{-1}(S^0),$$
we require that: $\pi$ is the projection on the second factor; each $E_{j+1}$ is the standard blowing-up centered at  $S^j$; each $S^j$ is a union of irreducible components of $\Sigma^j$; each $\Sigma^j$ is a smooth subset of the divisor $\DD^j$. Moreover, the maps $\pi_j|_{\Sigma^j}$ and $\pi_j|_{S^j}$ are etale over $K$. The set of irreducible components of $\DD^j$ is denoted by $\text{Comp}(\DD^j)$. The integer $h$ is called {\it height} of the blowing-up process and
$\left(\MM^h,\DD^h,\Sigma^h,\pi_h\right)$ {\it the top} of the process. The composed map $E_h$ is called the {\it total morphism} of the process.

\noindent More generally, we call {\it tree over $K$} a quadruplet
$(\MM,\DD,\Sigma,\pi)$ such that there exists a biholomorphism with the top of a blowing-up process over $K$
$$(\MM,\DD,\Sigma,\pi)\simeq \left(\MM^h,\DD^h,\Sigma^h,\pi_h\right).$$ For any $J\subset K$ we denote by $(\MM_J,\DD_J,\Sigma_J,\pi_J)$ the tree over $J$ obtained by simple restriction.

\subsubsection{The sheaves $\GG^n_Z$, $n\geq 0$.}
 From now on, we fix a marked tree $(\MM,\DD,\Sigma,\pi)$. In order to get through a technical difficulty, the tree is enhanced with a {\it cross}: let $E$ be the total morphism of $(\MM,\DD,\Sigma,\pi)$.
\begin{definition}[Cross]\label {Croix}
A cross on $\MM$ is the strict transform $Z=E^*Z_0$ of a single $Z_0=\fami{Z_1}$ or of a couple $Z_0=\fami{Z_1,Z_2}$ of germs of smooth transversal hypersurfaces along $0\times K$. We assume that any component $Z$ meets a unique irreducible component of $\DD$.
\end{definition}

\noindent We consider
$\text{Aut}^C(\mathcal{M},Z)$ the sheaf over $\DD$ of germs of automorphism defined in a neighborhood of $\DD$ such that
$$ \pi \circ \Phi=\pi, \quad \Phi|_\DD=\textup{Id},\quad \Phi|_Z=\textup{Id}\quad \text{et}\quad \Phi|_{\pi^{-1}(C)}=\textup{Id}.$$

\noindent Let $\OO_\MM$ be the sheaf over $\DD$ of germs of functions on $\MM$. Let $D$ be an irreducible component of the divisor $\DD$. We denote $I_D\subset \OO_\MM$ the ideal subsheaf of germs of function vanishing along $D$. Let us consider the subsheaf $\mathfrak{M}\subset \OO_\MM$ generated by the pull-back of the maximal ideal along $0\times K$. We consider a filtration of  $\OO_\MM$ defined by
$\mathfrak{M}^n_Z:=I_Z\cdot\mathfrak{M}^n,\  n\geq 1$. The sheaf ${\mathfrak{M}^{\petit{C}}}^n_Z$ is a sub-sheaf of $\mathfrak{M}^n_Z$ whose sections vanish along $\pi^{-1}(C)$.
We also have to consider the sheaf $\mathfrak{I}_Z\subset\OO_\MM$  defined by
$\mathfrak{I}_Z:=I_Z\cdot\prod_{D\in Comp(\DD)} I_{D}.$

 \medskip 
\noindent We call $n^{th}$ infinitesimal crossed tree the analytical space
$$\MM^{[n],Z}:=\big( \DD, \OO_\MM\left/{\mathfrak{M}^n_Z}\right. \big).$$
The neighborhood of order $0$ is 
$\MM^{[0],Z}:=\big( \DD, \OO_\MM\left/{\mathfrak{I}_Z}\right. \big).$
We also consider the following ringed spaces:
$
\MM^{\underline{n},Z} :=\big( \DD,\mathfrak{I}_Z\left/\mathfrak{I}_Z\mathfrak{M}^n_Z\right. \big)$ and $
\MM^{\underline{0},Z} :=\big( \DD,\mathfrak{I}_Z\left/\mathfrak{I}^2_Z\right. \big). 
$
 \begin{definition}We denote by $\textup{Aut}^C_n(\mathcal{M},Z)$ the subsheaf of $\textup{Aut}^C(\mathcal{M},Z)$ of germs that coincide with $\textup{Id}$ when restricted to 
the infinitesimal neighborhood of order $n$.
\end{definition}
\noindent The sequence of normal inclusions of groups  
$$ \ldots \triangleleft
\textup{Aut}^C_p(\mathcal{M},Z)\triangleleft\textup{Aut}^C_{p-1}(\mathcal{M},Z)\triangleleft\ldots\triangleleft\textup{Aut}^C_0(\mathcal{M},Z)
\triangleleft\textup{Aut}^C(\mathcal{M},Z)$$
points out the existence of a main part map for the filtration above
\begin{eqnarray}\label{superj0}
\textup{Aut}_n(\mathcal{M},Z)&\stackrel{\mathcal{J}_n}{\longmapsto}&\OO_{\MM^{[n],Z}},\\
\textup{Aut}_0(\mathcal{M},Z)&\stackrel{\mathcal{J}_0}{\longmapsto} &\OO_{\MM^{\underline{0},Z}}.
\end{eqnarray}
The description of the above morphims requires the study of local expression of sections of $\textup{Aut}^C_p(\mathcal{M},Z)$. We refer to $\cite{moi}$ where this description is performed in the non-parameter case. 
\begin{definition} We denote ${\GG^C}^n_Z$ the subsheaf of
  $\textup{Aut}^C_n(\mathcal{M},Z)$ kernel of the morphism $\mathcal{J}_n$.
\end{definition}
\noindent These sheaves are sheaves of Lie groups associated to some natural sheaves of Lie algebras we are about to define. 
\subsubsection{The tree gluing.}\label{treeglu}

Thanks to the sheaf $\text{Aut}^C(\MM,Z)$, we are going to introduce a process called
{\it gluing} on $\MM$. This construction will allow us to define a large class of trees with same divisor analytical type. These trees will inherit a canonical cross.

\medskip
\noindent Let us define a particular type of open covering of the divisor. Let $\UU=\left\{U_i\right\}_{i \in \mathbb{I}=\mathbb{I}_0\cup \mathbb{I}_1}$ be the covering of $\DD$ constituted of two kinds of open sets: if $i$ belongs to $\mathbb{I}_0$, $U_i$ is the trace on $\DD$ of a neighborhood of a unique singular locus conformally equivalent to a polydisc; if $i$ belongs to $\mathbb{I}_1$, $U_i$ is an irreducible component of $\DD$ deprived of the singular locus of $\DD$. Such a covering is called \emph{distinguished} when there is no $3$-intersection. Distinguished coverings contain Stein open sets having fundamental systems of Stein neighborhood. From now on, a covering
denoted by $\UU$ will always supposed to be distinguished. 

\medskip
\noindent Thanks to distinguished covering, we are able to glue the open sets of that covering by identifying points with respect to a $1$-cocycle in $\textup{Aut}^C(\mathcal{M},Z)$. Let 
$\coup{\phi_{ij}}$ be a $1$-cocycle in $\ZZ^1\left(\UU,\textup{Aut}^C(\mathcal{M},Z)\right)$. We define:
$$\MM[\phi_{ij}]= \bigcup_i \mathcal{U}_i\times\{i\}
/_{\big(\{x\}\times\{i\}\sim \{\phi_{ij}(x)\}\times\{j\}\big)}, $$
 where $\UU_i$ is a neigborhood of $U_i$ in $\MM$ such that $\phi_{ij}$ becomes an automorphism of $\MM$ along $U_i\cap U_j$. The latter automorphism is well defined since it coincides with $\text{Id}$ along the divisor. The manifold we get comes with an embedding 
\begin{equation}\label{inj}
\DD \hookrightarrow \MM[\phi_{ij}]
\end{equation} 
whose image is denoted by $\DD[\phi_{ij}]$ and $\MM[\phi_{ij}]$ is considered as a germ of neighborhood of $\DD[\phi_{ij}]$.

\begin{definition}
The manifold germ $\MM[\phi_{ij}]$ is called gluing of $\MM$ along $\UU$ by the cocycle $\coup{\phi_{ij}}$. 
\end{definition}
\noindent By construction of the sheaf $\textup{Aut}^C(\mathcal{M},Z)$, we have the following isomorphism of tree
$$\MM[\phi_{ij}]|_C\simeq \MM|_C$$

\noindent The gluing of a crossed tree comes naturally with a cross: it is the direct image of $Z$ by the quotient map for gluing relation. Such a tree and cross are respectively denoted by
\begin{eqnarray*}
\left(\MM[\phi_{ij}],\DD[\phi_{ij}],\Sigma[\phi_{ij}]\right),\text{ and } Z[\phi_{ij}].
\end{eqnarray*}

\medskip
\noindent We associate to any gluing the data of morphisms on infinitesimal neigborhood generalizing the embedding (\ref{inj}). The description of ${\GG^C}^n_Z$ sections reveals the following property
\begin{propr}\label{superinj} Let $n$ be an integer and $\NN=\MM[\phi_{ij}]$ be a gluing of $\MM$ by a cocycle in $\ZZ^1(\UU,{\GG^C}^n_Z)$. Then the canonical analytical embeddings 
\begin{eqnarray*}
\rho^{[n]}_{\NN}:&\MM^{[n],Z}\hookrightarrow \NN, \\
\rho^{\underline{n}}_{\NN}:&\MM^{\underline{n},Z}\hookrightarrow \NN
\end{eqnarray*}
have for respective images $\NN^{[n],Z[\phi_{ij}]}$ and $\NN^{\underline{n},Z[\phi_{ij}]}$.
\end{propr}

\subsubsection{The $\text{Glu}^C_n(\MM,Z,\UU)$ categories.}
Let $p$ be an integer. Let us consider the crossed tree built by a succession of gluings
\begin{equation}\label{var1}
\MM[\phi^1_{ij}][\phi^2_{ij}][\ldots][\phi^p_{ij}]
\end{equation}
where
\begin{itemize}
\item $\coup{\phi^1_{ij}}$ is a $1$-cocyle of $\GG^n_Z$;
\item for $k=2,\ldots,p$, $\coup{\phi^k_{ij}}\in\ZZ^1\left(\UU[\phi^1_{ij}]\ldots[\phi^{k-1}_{ij}],{\GG^C}^n_{Z[\phi^1_{ij}]\ldots[\phi^{k-1}_{ij}]}\right)$ with $$\GG^n_{Z[\phi^1_{ij}]\ldots[\phi^{k-1}_{ij}]}\subset \text{Aut}^C\left(\MM[\phi^1_{ij}]\ldots[\phi^{k-1}_{ij}],Z[\phi^1_{ij}]\ldots[\phi^{k-1}_{ij}]\right).$$ 
\end{itemize}
Following the (\ref{superinj}) lemma, we have canonical embeddings
\begin{eqnarray}
\rho^{[n]}_{\MM[\phi^1_{ij}][\phi^2_{ij}][\ldots][\phi^p_{ij}]}:&\MM^{[n],Z}\hookrightarrow \MM[\phi^1_{ij}][\phi^2_{ij}][\ldots][\phi^p_{ij}],\label{inj1} \\
\rho^{\underline{n}}_{\MM[\phi^1_{ij}][\phi^2_{ij}][\ldots][\phi^p_{ij}]}:&\MM^{\underline{n},Z}\hookrightarrow \MM[\phi^1_{ij}][\phi^2_{ij}][\ldots][\phi^p_{ij}]. \label{inj2}
\end{eqnarray}

\begin{definition}
The $\textup{Glu}^C_n(\MM,Z,\UU)$ category is the category whose objects are crossed trees built as (\ref{var1}) with the data of embeddings  (\ref{inj1}) and (\ref{inj2}). Arrows are biholomorphic germs that commute with embeddings. If $\MM$ and $\NN$ are isomorphic in $\textup{Glu}^C_n(\MM,Z,\UU)$, we denote
$$\MM\stackrel{\mathtt{G}_n}{\simeq} \NN.$$ 
\end{definition}

\subsection{Proof of the existence theorem.}\label{demons}

From now on, we assume the tree $\MM$ to be the support of an isoholonomical deformation $\FF$. We are going to define a cobordism notion in order to detect on any
element of $\text{Glu}^C_0(\MM,Z,\UU)$ the existence of an isoholonomical deformation linked to 
$\FF$ by isoholonomical deformations. 
\begin{definition}[Cross adapted to $\FF$]
\noindent Let $Z$ be a cross on $\MM$. $Z$ is said to be adapted to 
$\FF$ when each component $Z_i$ is either a separatrix of $\FF$ or is attached at a regular point of $\FF$. In the latter case, $Z_i$ will be transversal to the leaves of the isoholonomical deformation.
\end{definition}
\noindent Now, we consider the sheaf $\mathfrak{X}^C_{S,Z}$ over $\DD$ which is a subsheaf of the sheaf of holomorphic vector field. A section $X$ of $\mathfrak{X}^C_{S,Z}$ is supposed to vanish on $\pi^{-1}(C)$, to be tangent to $\DD$, to the separatrix and to the cross. We also assume that $X$ is \emph{vertical}:
$$D\pi(X)\equiv 0.$$
The sub-sheaf $\mathfrak{X}^C_{\FF,Z}\subset\mathfrak{X}^C_{S,Z}$ is the sheaf of vector fields tangent to the deformation $\FF$.
From now on,  $e^{tX}$ refers to the flow of the vector field $X$ at time $t$. One notices that, if $X$ is a section of $\mathfrak{I}_Z\mathfrak{X}^C_{S,Z}$ then, for all $t\in \mathbb{C}$, $e^{tX}$ exists as germ and defines a section of $\text{Aut}^C(\MM,Z)$.

\noindent 

\begin{definition}[Cobordism]
Let $\NN$ be in $\textup{Glu}^C_0(\MM,Z,\UU)$. 
$\NN$ is said to be
$\FF$-cobordant to $\MM$ if there exists a finite sequence of $1$-cocyles
$\coup{T^k_{ij}}_{k=1,\ldots,N}$ such that the two following conditions are verified:
\begin{enumerate}
\item for any $p=0,\ldots,N-1$, let $\mathfrak{X}^C_{\FF_p,Z_p}$ be the sheaf over $\DD[e^{T^1_{ij}}][\cdots][e^{T^{p}_{ij}}]$ of germs of vector field tangent to the following deformation and cross $$\FF_p=\FF[e^{T^1_{ij}}][\cdots][e^{T^{p}_{ij}}],\quad Z_p=Z[e^{T^1_{ij}}][\cdots][e^{T^{p}_{ij}}].$$ We assume $\coup{T^{p+1}_{ij}}$ is a $1$-cocycle with values in  $\mathfrak{X}^C_{\FF_p,Z_p}$.

\item $\NN\stackrel{\mathtt{G}_0}{\simeq}\MM[e^{T^1_{ij}}][\cdots][e^{T^N_{ij}}]$.
\end{enumerate}
We summarize this definition with the following notation:
$$
 \MM\fleche{\FF_1,Z_1}\MM_2\fleche{ \FF_2,Z_2}\cdots \fleche{ \FF_{N-1},Z_{N-1}}\MM_N\stackrel{\mathtt{G}_0}{\simeq} \NN.
$$
\end{definition}
\noindent Here, $\NN$ inherit a canonical isoholonomical deformation embedded in an isoholonomical deformation of $\FF$. For example if the cobordism is elementary, i.e $N=1$, the deformation of the ambient space defined by $t \rightarrow \MM_t=\MM[e^{(t)T_{ij}}],\ t\in\dbar$ carries an isoholonomic deformation.
\subsubsection{First step : the infinitesimal level.}\label{demo1}

Let $\FF_0$ be the isoholonomical deformation over $K$ such that $E^*\FF_0=\FF$. In view of the definition , the deformation $\FF_0$ is given by a $1$-form $\Omega_0$. Let us denote $F_0$ a reduced equation of the separatrix of $\FF_0$. Since the deformation $\FF$ is locally trivial, one can reproduce the computation done in \cite{moi} in the non-parameter case to prove the following lemma based upon the generalized curve hypothesis:
\begin{prop}\label{suiteexapara}
There exists an exact sequence of sheaves
$$0 \longrightarrow \mathfrak{M}^n_Z\mathfrak{X}^C_{\FF,Z} \longrightarrow
\mathfrak{M}^n_Z\mathfrak{X}^C_{S,Z}\xrightarrow{E^*\Omega_0(.)} {\mathfrak{M}_Z^\petit{C}}^n
\left(F_0\circ E\right)\longrightarrow 0 $$
where $(F_0\circ E)$ is the sub-sheaf of $\OO_\MM$ generated by the function $F_0\circ E$. 
\end{prop}
\noindent In order to establish an equivalent of the infinitesimal cobrdism result in \cite{moi} in the present parameter context, we only have to show
\begin{lemme}\label{cohonul}
$$H^1(\DD,{\mathfrak{M}^\petit{C}_{Z}}^n)=0.$$
\end{lemme}
\noindent \begin{demo}
Let us consider $\MM_{\WW}$ a neighborhood of $\DD$ and the associted fibration 
$$\Pi:\MM_{\mathcal{W}}\longmapsto \mathcal{K}=\pi\circ E$$
where $K\subset \mathcal{K}$ is a Stein open set. The spectral sequence of \cite{Godement} associated to the sheaf  $\mathfrak{M}^n_Z$ and the fibration $\Pi$ induce an exact sequence
$$ H^1\left(\mathcal{K},\Pi_*{\mathfrak{M}^\petit{C}_{Z}}^n\right)\rightarrow H^1\left(\MM_{\mathcal{W}},{\mathfrak{M}^\petit{C}_{Z}}^n\right)\rightarrow H^0\left(\mathcal{K},\mathcal{R}^1\Pi_*{\mathfrak{M}^\petit{C}_{Z}}^n\right)\rightarrow H^2\left(\mathcal{K},\pi_*{\mathfrak{M}^\petit{C}_{Z}}^n\right).$$
Since $\pi$ is proper and $\mathfrak{M}^\petit{C}_{Z}$ coherent, $\Pi_*{\mathfrak{M}^\petit{C}_{Z}}^n$ is a coherent sheaf \cite{grauerthans}. As $\mathcal{K}$ is Stein, each extremal term of the sequence vanishes \cite{grauertremmert}. Hence, we have
\begin{equation}\label{grauert0}
 H^1\left(\MM_\mathcal{W},{\mathfrak{M}^\petit{C}_{Z}}^n\right)\simeq H^0\left(K,\mathcal{R}^1\Pi_*{\mathfrak{M}^\petit{C}_{Z}}^n\right).
\end{equation}
The fiber of the derived sheaf satisfies: 
\begin{equation}\label{fibreenx}
\left(\mathcal{R}^1\Pi_*{\mathfrak{M}^\petit{C}_{Z}}^n\right)_x\simeq  H^1\left(\Pi^{-1}(x),{\mathfrak{M}^\petit{C}_{Z}}^n|_{\Pi^{-1}(x)} \right).
\end{equation}
Let us denote by $\mathfrak{M}_x\subset \coup{\OO_{\mathcal{K}}}_x$ the ideal of germ of function vanishing at $x$. Thanks to a distinguished covering of $\Pi^{-1}(x)$, one can see that
$$H^1\left(\Pi^{-1}(x),{\mathfrak{M}^\petit{C}_{Z}}^n \right)\otimes_{{\OO_{\mathcal{K}}}_x}{\OO_{\mathcal{K}}}_x\left/_{\mathfrak{M}_x}\right.\simeq H^1\left(\Pi^{-1}(x),{\mathfrak{M}^\petit{C}_{Z}}^n\otimes_{{\OO_{\mathcal{K}}}_x}{\OO_{\mathcal{K}}}_x\left/_{\mathfrak{M}_x}\right.\right).$$
Let $i_x$ be the embedding $\Pi^{-1}(x)\subset \MM$; a simple local computation ensures that $${\mathfrak{M}^\petit{C}_Z}^n\otimes_{{\OO_{\mathcal{K}}}_x}{\OO_{\mathcal{K}}}_x\left/_{\mathfrak{M}_x}\right.\simeq i_x^*{\mathfrak{M}^\petit{C}_{Z}}^n.$$ Now, in view of the non-parameter computation in \cite{moi}, the cohomology of $i_x^*\mathfrak{M}^\petit{C}_{Z}$ satisfies
$$H^1\left(\Pi^{-1}(x),i_x^*{\mathfrak{M}^\petit{C}_{Z}}^n\right)=0.$$
Hence, one gets
$$H^1\left(\Pi^{-1}(x),{\mathfrak{M}^\petit{C}_{Z}}^n \right)\otimes_{{\OO_{\mathcal{K}}}_x}{\OO_{\mathcal{K}}}_x\left/_{\mathfrak{M}_x}\right.= 0.$$
The Nakayama's lemma \cite{lang} and the relation (\ref{fibreenx}) ensure that the derived  sheaf $\mathcal{R}^1\Pi_*{\mathfrak{M}^\petit{C}_{Z}}^n$ is the trivial sheaf. Using (\ref{grauert0}), we have 
$$H^1\left(\MM_\mathcal{W},{\mathfrak{M}^\petit{C}_{Z}}^n\right)=0.$$
Finally, the lemma is obtained by taking the inductive limit on a familly of Stein neighborhood of $K$.
\end{demo}

\noindent The long exact sequence of cohomology associated to the short sequence (\ref{suiteexapara}) and to the covering $\UU$ give us the infinitesimal cobordism property
\begin{prop}[Infinitesimal cobordism]
The canonical map
$$H^1\left(\DD,\mathfrak{M}^n_Z\mathfrak{X}^C_{\FF,Z}\right) \longrightarrow H^1\left(\DD,\mathfrak{M}^n_Z\mathfrak{X}^C_{S,Z}\right)$$
is onto.
\end{prop}

\subsubsection{Second step : cobordism in $\text{Glu}^C_1(\MM,Z,\UU)$.}\label{demo2}
This section is devoted to prove the following proposition:
\begin{prop}\label{Glulg1}
Any element of $\textup{Glu}^C_1(\MM,Z,\UU)$ is elementary $\FF$-cobordant to $\MM$.
\end{prop}
\noindent The proof performed in \cite{moi} using an algorithm of Newton type can be reproduced in the $\textup{Glu}^C_1(\MM,Z,\UU)$ category. The last argument is a stability property for neighborhood of exceptional divisor. Hence, we have only to establish an equivalent of the latter result in our context.
\begin{thm}[Stability property]
For $n$ big enough, the image of the natural embedding
$$\textup{Glu}^C_n(\MM,Z,\UU)\hookrightarrow \textup{Glu}^C_1(\MM,Z,\UU)$$ contains trees isomorphic to $\MM$ in the category $\textup{Glu}^C_1(\MM,Z,\UU)$.
\end{thm}
\begin{proof} Let $p$ and $n$ be integers and $\NN=\MM[\phi_{ij}]$ with $\coup{\phi_{ij}}$ in $\ZZ^1\left(\UU,\textup{Aut}^C_{p}(\MM,Z)\right)$. For $p$ big enough, the image of the natural morphism $$H^1(\DD,\textup{Aut}_{p}(\MM,Z))\rightarrow H^1(\DD,\textup{Aut}_{n}(\MM,Z))$$ is trivial \cite{MatSal}. Hence, one can take a trivialisation of $\coup{\phi_{ij}}$
$$\phi_{ij}=\phi_i\circ \phi_{j}^{-1},\ \coup{\phi_i}\in \ZZ^0\left(\UU,\textup{Aut}^C_{n}(\MM,Z)\right).$$
By definition of the sheaf $\textup{Aut}^C_{p}(\MM,Z)$, by restriction over $C$, we find
$$\phi_i|_{\pi^{-1}(C)}=\phi_j|_{\pi^{-1}(C)}.$$
Hence, the familly $\coup{\phi_i}$ defines a global section over $C$ of $\textup{Aut}_{n}(\MM,Z)$. With the Hartogs's argument, this global section induces a germ of biholomorphism $\phi$ along $0\times K \subset \mathbb{C}^2\times \mathbb{C}^p$. This biholomorphism commutes with the projection and lets fixed each point of the cross. In some adapted coordinates $(x,y,t),\  (x,y)\in (\mathbb{C}^2,0),\ t\in K$, $\phi$ is written 
$$\phi(x,y,t)=(x,y,t)+H(x,y,t)\coup{\sum_{i,j\geq \nu}a_{ij}(t)x^iy^j),\sum_{i,j\geq \nu}b_{ij}(t)x^iy^j),0}$$
where $H$ is a reduced equation of the cross over $K$ and $a_{ij},b_{ij}$ holomorphic functions on $C$. Since $K$ is Stein, in view of \cite{Chabat} there exists functions $A_{ij}$ and $B_{ij}$ holomorphic on $K$ such that 
$$A_{ij}|_C=a_{ij}\quad\text{et}\quad B_{ij}|_C=b_{ij}.$$
Then, the map
\begin{equation}\label{prolongement}
\Phi(x,y,t)= (x,y,t)+H(x,y,t)\coup{\sum_{i,j}A_{ij}(t)x^iy^j,\sum_{i,j}B_{ij}(t)x^iy^j,0}
\end{equation}
is a germ of automorphism along $0\times \mathcal{K}\subset \mathbb{C}^2\times \mathcal{K}$ wihch extends $\phi$, commutes with the projection and fixes the cross. Morevover, if one chooses $A_{ij}=0$ et $B_{ij}=0$ as soon as $a_{ij}=0$ and $b_{ij}=0$, then the tangency order to the idenity of the extension is the same as $\phi$. Hence, for $p$ big enough, the biholomorphism (\ref{prolongement}) can be lifted up in a global section of $\textup{Aut}_n(\MM,Z)$ over $K$ with $\Phi|_C=\phi$. Hence, the $0$-cocycle $\coup{\psi_i}=\coup{\phi_i\circ\Phi^{-1}|U_i}$ is a trivialisation of $\coup{\phi_{ij}}$ with values in $\textup{Aut}^C_n(\MM,Z)$. 

\noindent Now, since $K$ is Stein, the curve $Z_0$ defining the cross $Z=E^*Z_0$ can be straightened along $0\times K$ and in some coordinates $(x,y,t),\ (x,y)\in \mathbb{C}^2, t\in K$, $Z_0$ admits $xy=0$ for equation. The latter coordinates induce two canonical systems of coordinates along the components of $Z=Z_1\cup Z_2$. The total morphism $E$ is now written 
$E(x_1,y_1,t)=(x_1,y_1x_1^{N_1},t)$ and
$E(x_2,y_2,t)=(x_2y_2^{N_2},y_2,t)$. Let us denote by $\psi_1$ and $\psi_2$ the components of $\coup{\psi_i}$ defined on the open set of the covering  which contains the components of $Z$. Since each component of $Z_0$ is smooth, the automorphisms $\psi_1$ and $\psi_2$ can be written
\begin{eqnarray*}
\psi_1(x_1,y_1,t)&=&\left(x_1+x_1^ny_1U_1(x_1,y_1,t),y_1+x_1^ny_1V_1(x_1,y_1,t) ,t\right) \\
\psi_2(x_2,y_2,t)&=&\left(x_2+y_2^nx_2U_2(x_2,y_2,t),y_2+y_2^nx_2V_2(x_2,y_2,t) ,t\right).
\end{eqnarray*}
Let $\psi$ be the germ of biholomorphism along $0\times K$ defined by
$$\psi(x,y,t)=\coup{x(1+y^nU_2(0,y,t)),y(1+x^nV_1(x,0,t)),t}.$$ 
For $n$ big enough, $\Psi$ can be lifted up on $\MM$ in an automorphism $\Psi$ that fixes each point of $\DD$ and $Z$. Now, if one takes a closed look to the expression of $\psi$, one can verify that, for any point $x$ on $\DD$ and any component $\psi_{i(x)}$ of $\coup{\psi_i}$ defined near $x$,
$$(\mathfrak{J}_1)_x(\psi_{i(x)})=(\mathfrak{J}_1)_x(\Psi).$$
Hence, the $0$-cocycle $\coup{\psi_i\circ \Psi^{-1}|_{U_i}}$ is a trivialisation of $\coup{\phi_{ij}}$ in $\GG^1_Z$. Since, the tree $\NN$ is the gluing $\MM[\phi_{ij}]$, $\NN$ is isomorphic to $\MM$ in the category $\textup{Glu}^C_1(\MM,Z,\UU)$.
\end{proof}

\subsubsection{Third step: cobordism in $\text{Glu}^C_0(\MM,Z,\UU)$.}\label{demo3}
The cobordism in $\text{Glu}^C_0(\MM,Z,\UU)$ is related to the following result: 
\begin{prop}\label{cobordisme0}
Any element in $\textup{Glu}^C_0(\MM,Z,\UU)$ is $\FF$-cobordant to $\MM$.
\end{prop}
\noindent The proof done in \cite{moi} for the category $\textup{Glu}_0(\MM,Z,\UU)$ can be repeated here without any change. The main tools are an induction on the height of the trees and the cobordism result for $\text{Glu}^C_1(\MM,Z,\UU)$.

\subsubsection{Fourth step: preparation of a cocyle.}
Let $(\MM',\DD',\Sigma',\pi')$ be a tree over $K$ topologically equivalent to $(\MM,\DD,\Sigma,\pi)$. We suppose that over $C$ the trees $\MM'|_C$ and $\MM|_C$ are analytically equivalent. 
\begin{prop}\label{supergoldo} There exists an isoholonomic deformation $\FF'$ on $\MM'$ such that $\FF'|_C$ and $\FF|_C$ are analytically equivalent and $\FF$ and $\FF'$ are embedded in an isoholonomic deformation over $K\times \dbar$ which is trivial along any $\{t\}\times \dbar,\ t\in C$.
\end{prop}
\noindent In order to prove the latter proposition, we want to prepare a $1$-cocycle such that the tree $\MM'$ becomes a gluing of $\MM$ in a category $\textup{Glu}_0(\MM,Z,\UU)$. This cocycle must well behave with respect to the condition $\MM'|_C\simeq\MM|_C$. 

\medskip
\noindent In view of \cite{Seguy}, one can find an isoholonomic deformation $\tilde{\FF}$ over $\dbar\times K$ such that:
$\tilde{\FF}|_{-1\times K}$ is equal to $\FF$; $\tilde{\FF}|_{1\times K}$ is a deformation taken on a tree $(\tilde{\MM},\tilde{\DD},\tilde{\Sigma},\tilde{\pi})$ with $\DD'$ and $\tilde{\DD}$ analytically equivalent; the deformations $\tilde{\FF}|_{1\times C}$ and $\FF|_C$ are\hyphenation{ana-ly-ti-cally} analytically equivalent. Let us denote by $\theta$ a conjugacy between $\DD'$ and $\tilde{\DD}$ and $\phi$ a conjugacy between $\MM'|_C$ and $\tilde{\MM}|_C$. We point out the fact that, generally, $\theta$ cannot be extend to the trees. Hence, we are far from getting a deformation on $\MM'$. However, since $K$ is Stein, the tubular neighborhood of any irreducible component of $\DD$ is a trivial deformation over $K$. Therefore, for such any component $D$, there exists a biholomorphism $\Theta_D$ from a tubular neighborhood $T(D)$ of $D$ to a tubular neighborhood of $\theta(D)$ extending $\theta|_D$ such that $\Theta_D(\tilde{\DD})=\DD'$. The automorphism of $T(D)|_C$ defined by $$\Theta_D|_C\circ \phi^{-1}|_{T(D)}$$ is well-defined on a neigborhood of $D|_C$, lets invariant each component of $\DD|_C$ and commutes with the projection.
\begin{lemme}\label{coborauto}
There exists an automorphism $\Phi_D$ of a neighborhood of $D$ extending $\Theta_D|_C\circ \phi^{-1}|_{T(D)}$ over $K$, which lets invariant $\DD$ and commutes with $\pi$.
\end{lemme}
\noindent \begin{demo} Let us denote by $h_D$ the automorphism $\Theta_D|_C\circ \phi^{-1}|_{T(D)}$. Since $K$ is Stein, there exists a system of coordinates $(x,t,s),\  x\in \mathbb{P}^1,\ t\in (\mathbb{C},0),\ s\in K$ in a neighborhood of $D$ such that: $\{x=0\}$ is a local equation of $D$; the components of $\DD$ transversal to $D$ have equations of the form $\fami{ t=f_i(s)  , t=\infty }_{i=1,\ldots,N}$; the fibration $\pi$ is $\pi:(x,t,s)\mapsto s$. In view of all its properties, $h_D$ is written

$$(x,t,s)\mapsto \coup{xA(x,t,s),t+xB(x,t,s)\prod_{i=1,\ldots,N}(t-f_i(s)),s},\quad s\in C$$
where $A$ and $B$ are holomorphic functions. Since $h_D$ is a germ of automorphism, we have $A(0,t,s)\neq 0$. Moreover as $h_D$ is global and extendable along $\{t=\infty \}$, if $A$ and $B$ are written
$$A(x,t,s)=\sum_{ij}a_{ij}(s)x^it^j,\quad B(x,t,s)=\sum_{ij}b_{ij}(s)x^it^j,\quad s\in C$$
the functions $a_{ij}$ and $b_{ij}$ vanish as soon as $ip-j<0$, $p$ refering to the self-intersection of $D$. Particulary, we find $A(0,t,s)=a_{00}(s)$ for any $s,t$. Since $K$ is Stein and compact, $a_{00}$ can be extended in a non-vanishing holomorphic function on $K$ and any other function $a_{ij}$ or $b_{ij}$ can be extended too. If one carefully chooses to extend with the zero function as soon as $a_{ij}$ or $b_{ij}$ is the zero function, one gets an extension on $h_D$ satisfying all the checked properties.
\end{demo}
\noindent Now for any component $D$, let us consider
$\Lambda_D=\phi_D^{-1}\circ \Theta_D.$
If we glue the familly of tubular neighborhoods with respect to the familly of automorphisms $\Lambda_{DD'}=\Lambda^{-1}_D\circ \Lambda_{D'}$, we find:
\begin{eqnarray}
\MM'\simeq \left.\coprod_{D\in \text{Comp}(\tilde{\DD})} T(D)\right/\coup{x\sim \Lambda_{DD'}(x)}_{(D,D')\in \text{Comp}(\tilde{\DD})^{\check{2}}},\label{collageenorme} \\
\Lambda_{DD'}|_C=\text{Id}.\label{coherence}
\end{eqnarray}
From now on, we make two operations on the familly $\coup{\Lambda_{DD'}}$ which leads us to a $1$-cocyle taking its values in ${\GG^C}^0_Z$. Using an analogous with parameter of the lemma $(3.3)$ in \cite{moi} , one can first suppose $\Lambda_{DD'}$ to be tangent to the identity along the singular locus $D\cap D'$. Then, by taking a distinguished covering $\UU$ of $\DD$ finer than the tubular neighborhood as in $(3.2)$ of \cite{moi}, one can build an element $\coup{\phi_{ij}}$ related to $\Lambda_{DD'}$, which belongs to $\ZZ^1(\UU,{\GG^C}^0_Z)$ for a cross $Z$ well chosen. Moreover, in this construction, one keeps the property $\MM'\simeq \tilde{\MM}[\phi_{ij}]$. Hence, the tree $\MM'$ is conjugated to an element of $\textup{Glu}_0(\MM,Z,\UU)$. Therefore, the proposition is a consequence of (\ref{cobordisme0}).

\subsubsection{Fifth step: a finite determinacy argument.}
If $E$ is the total morphism of a blowing-up process, let us denote by $\text{Att}(E,S)$ the locus of intersection between the exceptional divisor and the strict transform of $S$ by $E$. Let us prove now the theorem (\ref{exis}). Let the tree $(\MM^1,\DD^1,\Sigma^1,\pi^1)$ be the top of the desingularization process of the equisingular family $S^1$. For $n$, let us denote by $(\MM_n^1,\DD_n^1,\Sigma_n^1,\pi_n^1)$ the top of the process given by the following diagram
\begin{equation}
\begin{array}{ccccccccc}
\MM_n^1 & \stackrel{E^n}{\rightarrow} & \ldots &\MM^1_{j} & \stackrel{E^j}{\rightarrow} &  & \ldots & \stackrel{E^1}{\rightarrow} & \MM^{1}  \\
\bigcup &  &    &\bigcup &  &
 &  &  & \bigcup  \\
\text{Att}(E_n,S^1) & \rightarrow & \ldots &\text{Att}(E_j,S^1) & \rightarrow &  & \ldots & \rightarrow & S^1\\
\bigcup &  &  &  \bigcup &  &  & & & \bigcup   \\
\text{Att}(E_n,S^1) & \rightarrow & \ldots &\text{Att}(E_j,S^1) & \rightarrow &  & \ldots & \rightarrow & S^1

\end{array}\end{equation}
where $E_j$ refers to the total morphism of $\MM_{j-1}^1$. For any integer $n$, the proposition (\ref{supergoldo}) ensures the existence of an isoholonomic deformation $\FF_n^1$ on $\MM_n^1$ such that $\FF_n^1|_C$ and $\FF^0|_C$ are analytically equivalent. For $n$ big enough, one can suppose that each component of $\DD_n^1$ meets at most one irreducible component of $E_n^*S^1$. Moreover, $S^1$ and $S^0$ are topologically equivalent as well as $\FF_n^1$ and $E_n^*\FF^0$. Hence, each component of $\DD_n^1$, which meets a component of $E_n^*S^1$, meets exactly one component of the separatrix of $\FF_n^1$. In view of \cite{mather}, a property of finite determinacy ensures that for $n$ big enough, the separatrix of $\FF_n^1$ and $E_n^*S^1$ are analytically equivalent. Therefore, $\FF_n^1$ can be pulled down in an isoholonomical deformation $\FF^1$ over $ K$ satisfiying $\text{Sep}(\FF^1)=S^1$ and $\FF^1|_C\simeq \FF^0|_C$. Moreover, by construction, $\FF^0$ and $\FF^1$ are embedded isoholonomical deformations over $K\times \dbar$ satisfying the checked properties.

\section{Quasi-homogeneous foliations.}\label{part2}

\subsection{Desingularization of quasi-homogeneous foliations.}
\noindent We are interested first in the desingularization of quasi-homogeneous foliation of generalized curve type.

\noindent Let $\Omega$ be an open set in $\mathbb{C}^2$, $p$ a point in $\Omega$ and $S$ a germ of smooth curve. For any proper morphism $E:X\rightarrow \Omega$, we call the \emph{strict transform of $S$ by $E$} the closure in $X$ of the analytical set $E^{-1}(S\backslash 0)$. The intersection of the strict transform of $S$ and the divisor $E^{-1}(0)$ is called the attaching locus of $S$ with the respect to $E$.

\noindent We define a special morphism with a finite number of successive blowing-up centered at point $\mathfrak{E}^n_d(p,S):(\MM^n_d,\DD^n_d)\rightarrow (\Omega,p)$ by  $$\mathfrak{E}^n_d(p,S)=E_1\circ \cdots\circ E_n$$ where $E_i$ is the blowing-up centered at the attaching locus of $S$ with respect to $E_1\circ \cdots \circ E_{i-1}$. The exceptional divisor $\DD^n_d=\mathfrak{E}^n_d(p,S)^{-1}(0)$ is a chain of $n$ irreducible components $\{D_i\}_{1\leq i\leq n}$ such that $$(E_1\circ\ldots E_i)^{-1}(p)=D_1\cup\ldots \cup D_i.$$ The component $D_i$ is called the $i$-component of $\mathfrak{E}^n_d(p,S)$

\smallskip
\noindent Let $f$ be a reduced quasi-homogeneous polynomial function defined in some coordinates $(u,v)$
$$ f=\left\{\sum_{\alpha i + \beta j = \gamma} a_{ij}u^iv^j=0\right\},\quad \alpha,\beta,\gamma \in \mathbb{N}^*,\quad \text{pgcd}(\alpha,\beta) =1,\ \alpha<\beta  .$$
Let us write the Euclide algorithm for the couple $(\alpha,\beta)$:
\begin{equation}\label{Euclid} r_0=\beta,\ r_1=\alpha \quad
\left\{\begin{array}{c}
 r_0=q_1r_1 + r_2 \\
\cdots \\
r_i=q_{i+1}r_{i+1}+r_{i+2} \\
\cdots\\
r_N=q_{N+1}r_{N+1}+0\\
\end{array}\right. .
\end{equation}
We are now able to give a precise description of the morphism of desingularization of $f$.
\begin{prop} The morphism of desingularization of $f^{-1}(0)$ is the composition
\begin{equation}\label{morphism}
\mathfrak{E}^{q_1}_d(p_1,S_1)\circ\mathfrak{E}^{q_2}_d(p_2,S_2)\circ\cdots\circ\mathfrak{E}^{q_{N+1}}_d(p_{N+1},S_{N+1})
\end{equation}
where 
\begin{enumerate}
\item $p_1=0,\  S_1=\{v=0\}$.
\item $p_{i+1}$ is the intersection point of the $(q_{i})$-component and the $(q_{i}-1)$-component of $\mathfrak{E}^{q_{i}}_d(p_{i},S_{i})$.
\item $S_{i+1}$ is the germ of smooth curve defined by the $q_{i-1}$-component of $\mathfrak{E}^{q_{i}}_d(p_{i},S_{i})$.
\end{enumerate}
\end{prop}
\begin{proof}
The proof is an induction on the length of the Euclid algorithm. If the length is $1$, since $\alpha$ and $\beta$ are relatively prime, the algorithm is reduced to
$$r_0=\beta,\ r_1=\alpha=1 \quad r_0=r_0\times r_1. $$
In the coordinates $(u,v)$, the morphism $\mathfrak{E}^{r_0}_d(0,\{v=0\})$ is locally written
$$\mathfrak{E}^{r_0}_d(0,\{v=0\})(u_1,v_1)=(u_1,v_1u_1^{r_0}).$$
Hence, the pull-back function  $\mathfrak{E}^{r_0}_d(0,\{v=0\})^*f$ is expressed by
$$\mathfrak{E}^{r_0}_d(0,\{v=0\})^*f(u_1,v_1)=\sum_{ i + r_0 j = \gamma} a_{ij}u_1^{i+r_0j}v_1^j=u_1^{\gamma}\sum_{ i + r_0 j = \gamma} a_{ij}v_1^j$$
The curve $\sum_{ i + r_0 j = \gamma} a_{ij}v_1^j$ is the union of a finite number of smooth curves transversal to the exceptional divisor. One can verify that in any other canonical coordinates for the blowing-up morphism $\mathfrak{E}^{r_0}_d(0,\{v=0\})$, the function $\mathfrak{E}^{r_0}_d(0,\{v=0\})^*f$ does not vanish along any curve transversal to the exceptional divisor except maybe on the first component: the latter case only occurs when the axe $\{u=0\}$ in the considered coordinates is an irreducible component of $f^{-1}(0)$. In any case, the curve $f^{-1}(0)$ is desingularized by $\mathfrak{E}^{r_0}_d(0,\{v=0\})$.
\noindent If the Euclid algorithm is of length $n$, using the notation of (\ref{Euclid}), one considers the morphism 
$\mathfrak{E}^{q_1}_d(p_1,S_1)$. As above, the pull-back equation is written 
\begin{eqnarray*}f_1(u_1,v_1)=\mathfrak{E}^{q_1}_d(0,\{v=0\})^*f(u_1,v_1)&=&\sum_{\alpha i + \beta j = \gamma} a_{ij}u_1^{i+q_1j}v_1^j\\
&=&\sum_{ \alpha(i + q_1 j)+r_2j = \gamma} a_{ij}u_1^{i+q_1j}v_1^j.
\end{eqnarray*}
To conclude, one observes that the last equation defines a quasi-homogeneous curve of weight $(\alpha,r_2)$. Since the Euclid algorithm for the couple $(\alpha,r_2)$ is of length $n-1$, the induction hypothesis ensures that the desingularization of the latter curve admits a description of  type (\ref{morphism}). As the desingularization of $f^{-1}(0)$ is the successive composition of $\mathfrak{E}^{q_1}_d(0,\{v=0\})$ and of the desingularization of $f_1^{-1}(0)$, the proposition is proved.
\end{proof}
\noindent Now, if $\FF$ is a quasi-homogeneous curve of generalized curve type, its desingularization admits the same description as (\ref{morphism}). In particular, $\FF$ has a dual tree of the form depicted on Figure \ref{arbre}~:
\begin{figure}[!ht]
\center{\includegraphics{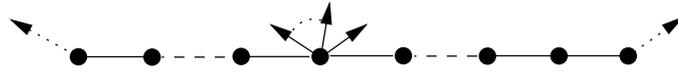}}
\caption{Dual tree of a quasi-homogeneous foliation}\label{arbre}
\end{figure}
the non-extremal vertex with some arrows corresponds to the central component. The extremal components carry an arrow if and only if, in the $(u,v)$ coordinates, at least one of the axes is an irreducible component of the separatrix.

\subsection{Automorphism of reduced singularity.}
In order to prove that there is no \emph{transversal obstruction} for quasi-homogeneous foliations with same holonomy to be conjugated and in view of the cohomological interpretation (\ref{intco}), we first recall some classical facts about automorphisms of reduced singularities. 

\bigskip
\noindent Let $\FF$ be a germ of reduced singularity with two non-vanishing eigen-values. In some coordinates \cite{MM}, the foliation $\FF$ is given by an holomorphic $1$-form $\omega$ where
\begin{equation}\label{normalform}
\omega = x(1+A(x,y))\dd + y(\lambda +B(x,y))\dd y
\end{equation}
with $A(0,0)=B(0,0)=0$. Hence, the axes $\{x=0\}$ and $\{y=0\}$ are both separatrix. In fact, these are the sole separatrix. Let us denote by $S$ and $S'$ the respective axes.
For any open set $U$ in $S$, the notation $\text{Aut}(\FF,U)$ refers to the group of local automorphism germs $\phi$ along $U$, which let globally invariant the foliation and satisfy
$$\phi|_U=\textup{Id}.$$ Precisely, a germ of automorphism $\phi$ is in $\text{Aut}(\FF,U)$ if and only if:
$$\phi^*\omega\wedge \omega=0,\quad \phi|_U=\textup{Id}.$$
Let us denote by $\text{Fix}(\FF,U)$ the sub-group of $\text{Aut}(\FF,U)$, which lets invariant each local leaf. In \cite{iso}, D. Cerveau and R. Meziani call $\text{Aut}(\FF,U)$ the \emph{isotropy group} of the singularity. They give a description of all elements in $\text{Fix}(\FF,U)$ when $U$ is a neighborhood of $0$ in the separatrix $S$. To be more specific, any element of $\text{Fix}(\FF,U)$ can be written $$(x,y)\mapsto e^{(\tau(x,y))X}$$ Here, $X$ is a germ of tangent vector field and $\tau(x,y)$ a germ of holomorphic function vanishing along $\{x=0\}$. The notation $e^{(t)X}$ refers to the flow of $X$ at time $t$. This kind of description persists if one considers a punctured neighborhood of $0$ in $S$ or even any corona around $0$.

\bigskip
\noindent Let $T$ be a germ of curve transversal to $S$ and $\text{Hol}_T\in \textup{Diff}(T,T\cap S)$ be the holonomy automorphism of $\FF$ computed on $T$. Let $U$ be either a neighborhood of $0$ in $S$ or a corona around $0$ in $S$ such that $T\cap S\in U$. Let $\phi$ be any element of $\text{Aut}(\FF,U)$. In a small neigborhood of $T\cap S$ , the foliation $\FF$ can be straightened: precisely, there exists local coordinates $(u,v)$ such that $S=\{v=0\}$, $T=\{u=0\}$ and $\FF=\{v=\textup{cst}\}$. The curve $\phi(T)$ is transversal to $S$ and meets $S$ at $S\cap T$. Hence, there exists a germ of biholomorphism $\rho_{\phi}$ such that
$$\rho_{\phi}:(0,v)\in (T,S\cap T)\mapsto (\alpha(v),v) \in (\phi(T),S\cap T).$$
where $\alpha$ is a germ of holomorphic function with $\alpha(0)=0$ and $\alpha'(0)\neq 0$.
By definition, the holonomy map $\text{Hol}_{\phi(T)}$ computed as a germ in $\textup{Diff}(\phi(T),T\cap S)$ satisfies the following relation
$$ \text{Hol}_T=\rho_{\phi}^{-1}\circ \text{Hol}_{\phi(T)}\circ \rho_{\phi}$$
Let us consider a path $\gamma$ in the leaf passing through $(0,v)$ from the point $(0,v)$ to the point $(0,\text{Hol}_T(v))$ making exactly one turn around $S'$. Since $\phi$ is the identity when restricted to $S$, for $v$ small enough, $\phi(\gamma)$ is a path in the leaf from  $(0,\phi(v))$ to $(0,\phi(\text{Hol}_T(v)))$ making exactly one turn around $S'$ with the same orientation as $\gamma$. Hence, by definition of the holonomy,
$$ \text{Hol}_{\phi(T)}(\phi(v))=\phi(\text{Hol}_T(v)).$$
With the above relations, we find
$$(\rho_{\phi}^{-1}\phi)\circ\text{Hol}_T=\text{Hol}_T\circ(\rho_{\phi}^{-1}\phi)$$
Hence, we build a morphism defined by
$$\phi\in\textup{Aut}(\FF,U)\longmapsto \left[\rho_{\phi}^{-1}\phi\right]\in \textup{Cent}(\text{Hol}_T)/<\text{Hol}_T>.$$
It is a morphism of groups with values in the quotient of the centralisator $\textup{Cent}(\text{Hol}_T)$ by the abelian sub-group generated $\text{Hol}_T$. Moreover, we have the following result:

\begin{lemme}\label{Fix} If $U$ is a small enough neighborhood of the singularity or a corona around it then
the following sequence
\begin{equation}\label{morphisme}
0\longrightarrow \textup{Fix}(\FF,U)\longrightarrow\textup{Aut}(\FF,U)\longrightarrow \textup{Cent}(\textup{Hol}_T)\left/_{<\textup{Hol}_T>}\right.\longrightarrow 0
\end{equation}
is exact.
\end{lemme}
\noindent In the formal context, since one can define the holonomy of a transversally formal foliation along an irreducible component of the divisor, the lemma (\ref{Fix}) can be reproduced by using transversally formal vocabulary. 

\subsection{Proof of the convergent statement.}\label{findelapreuve}
One can clearly suppose that $\FF_0$ and $\FF_1$ are two quasi-homogeneous generalized curves with $$\textup{Sep}(\FF_0)=\textup{Sep}(\FF_1)=S.$$ As each foliation is of generalized curve type, their desingularizations are both equal to the desingularization of $S$ \cite{camacho}. Let us denote by $E:(\MM,\DD)\rightarrow (\mathbb{C}^2,0)$ the morphism of desingularization where $\DD$ refers to the exceptional divisor $E^{-1}(0)$. Since the projective holonomy of $E^*\FF_0$ and $E^*\FF_1$ over the central component $D_0$ are conjugated, the restricted foliation $E^*\FF_0$ and $E^*\FF_1$ on a tubular neighborhood of $D_0$ are conjugated: in order to prove the latter fact, let us consider a curve $T_0$ transversal to $D_0$, on which one computes the projective holonomy representations. Let the set $\fami{s_1,\ldots,s_n}$ refer to the singularities of $E^*\FF_0$ along $D_0$. In view of the hypothesis, there exists a germ of automorphism $h_0:T_0\rightarrow T_0$ so that for any $[\gamma]\in \Pi_1(D_0\backslash \{s_1, \ldots, s_n\})$ and $x\in T_0$, one has
\begin{equation}\label{holconj}
h_0([\gamma]_{\FF_0} x)=[\gamma]_{\FF_1}h_0(x)
\end{equation}
where $[\gamma]_{\FF_0}$ and $[\gamma]_{\FF_1}$ refer to the image of the path $\gamma$ through the respective projective holonomy representations. In view of the desingularization process of $\FF_0$, there exists a fibration $\pi$ over $D_0$, for which the irreducible components of the strict transform of $S$, which are attached to $D_0$, and the two components of $\DD$ transversal to $D_0$ are some fibers. Let us call them the \emph{special} fibers. One can choose the fibration $\pi$ such that any fiber different from the special fibers, is transversal to the leaves of the foliations $E^*\FF_0$ and $E^*\FF_1$, at least on a little neighborhood of $D_0$. Let us choose for $T_0$ any fiber of $\pi$ different from the special fibers.
For any point $x$ in a neighborhood of $D_0$ deprived of the special fibers, one can consider a path $\gamma(t),\ t\in[0,1]$, which links $x$ to some point of $T_0$ in the leaf. The point $H_0(x)$ is now defined as the extremity of the lifting path in the leaf passing by $\gamma(1)$ with respect to the fibration $\pi$. The relation (\ref{holconj}) ensures that this construction does not depend on the path $\gamma$; hence, $x\rightarrow H_0(x)$ is well defined. Moreover, by construction $H_0$ is bounded near the special fibers. Hence, $H_0$ can be holomorphically extended on a neighborhood of $D_0$. One can check that $H_0$ sends any local leaf of $E^*\FF_0$ on a local leaf of $E^*\FF_1$. We can observe that, by construction, the restriction of $H_0$ on the component $D_0$ is the identity.

\noindent  Let us denote by $D_1$ an irreducible component of the divisor meeting $D_0$ at the point $s_{01}$. Since $D_1$ has only two singular points, the holonomy representations are morphisms of the form  
$$k\in\mathbb{Z}=\Pi_1\big( D_1 \backslash\text{Sing}(\DD)\big)\rightarrow [\gamma_i]^k \in \text{Diff}(\mathbb{C},0)$$
where $\gamma_i$ is the holonomy of a path in $D_1$ making one turn around $D_0$ for the foliation $E^*\FF_i$. As the two foliations are analytically equivalent near $s_{01}$, the holonomy maps $[\gamma_0]$ and $[\gamma_1]$ are conjugated by an interior automorphism of $\text{Diff}(\mathbb{C},0)$. Hence, the whole projective holonomy representations over $D_1$ of $E^*\FF_0$ and $E^*\FF_1$ are conjugated. In view of the geometry of the desingularization process, one can repeat the argument for any component of the divisor. Then thanks to the same construction as before, one can extend any conjugacy of the holonomy over $D$ on a tubular neighborhood of $D$ denoted by $T(D)$. Hence, we get a germ of biholomorphism $H_D$ along each $D$ such that 
$$H_D^*E^*\FF_0|_{T(D)}=E^*\FF_1|_{T(D)}\quad\text{and}\quad H_D|_D=\text{Id}$$
Of course, there is no chance for the family $\fami{H_D}_{D\in \text{Comp}(\DD)}$ to induce a global biholomorphism, i.e to verify the condition
$$H_D=H_{D'},\ \text{on a neighborhood of } D\cap D'.$$
However, we are going to introduce a finer covering than $\fami{T(D)}_{D \in \text{Comp}(\DD)}$, which allows us to \emph{twist} the family $\fami{H_D}_{D\in \text{Comp}(\DD)}$ and to build a new family $\fami{\tilde{H}_i}_{i\in\mathbb{I}}$ satisfying 
$$\tilde{H}_i\circ \tilde{H}^{-1}_j\text{ acts in the local leaf of } E^*\FF_0.$$
Precisely, let us consider a covering $\fami{U_i}_{i\in \mathbb{I}=\mathbb{I}_0\cup\mathbb{I}_1}$ of $\DD$ defined by 
$$\left\{ \begin{array}{ll}
i \in \mathbb{I}_0= \textup{Comp}(\DD), & U_D=D\backslash \textup{Sing}(E^*\FF_0) \\
i \in \mathbb{I}_1=\textup{Comp}(\DD)^{2}, & U_{(D,D')}= T(D)\cap T(D')\cap \DD
\end{array}\right.$$
In view of the form of the dual graph of $\FF_0$ or $\FF_1$, we use the following clear notation for the components of $\DD$ 
$$\textup{Comp}(\DD)=\fami{D_{-m},D_{-m-1},\ldots,D_{-1},D_0,D_1,\ldots,D_{p-1}, D_p}$$
where $D_0$ refers to the central component. We consider the filtration $\mathbb{I}_n$ of $\mathbb{I}$ defined by 
$$\mathbb{I}_n=\fami{D_{-n},\ldots,D_n}\bigcup \fami{D_{-n},\ldots,D_n}^2$$
Using the special geometry of the dual tree of quasi-homogenous foliation, we establish the following lemma: 
\begin{lemme} For any integer $n$, there exists a family $\fami{\phi_i}_{i\in I_n}$ such that 
\begin{itemize}
\item for all $D\in \mathbb{I}_n$, $\phi_D$ belongs to $\textup{Aut}(E^*\FF_0,U_D)$ 
\item for all $(D,D')\in\mathbb{I}_n$, $\phi_{(D,D')}$ belongs to $\textup{Aut}(E^*\FF_0,U_{(D,D')})$ 
\end{itemize}
and such that for all $\coup{D,(D,D')}\in \mathbb{I}_0\times \mathbb{I}_1$ the two maps
\begin{eqnarray*}
\phi_D^{-1}\circ H_D^{-1}\circ H_{D'}\circ\phi_{(D,D')} \\
\text{and}\quad \phi_{D'}^{-1}\circ\phi_{(D,D')}
\end{eqnarray*}
act in the local leaf of $E^*\FF_0$.
\end{lemme}
\noindent \begin{demo}
The proof is an induction on the integer $n$. For $n=0$, the lemma is trivial: since the condition is empty, on can choose $\Phi_{D_0}=\text{Id}$. Let us suppose the result true for $n$. The automorphism $\phi_{D_n}^{-1}\circ H_{D_n}^{-1}\circ H_{D_{n+1}}$ is an automorphism of the foliation defined in a neighborhood of $U_{D_n}\cap U_{(D_n,D_{n+1})}$. By construction, the foliation $E^*\FF_0$ restricted to a neighborhood of $U_{(D_n,D_{n+1})}$ has exactly one reduced isolated singularity. Hence, one can consider $\text{Hol}_n$ its holonomy map computed on any curve transversal to $D_n$ and attached to a point in $U_{(D_n,D_{n+1})}$. In view of the lemma (\ref{Fix}), there exists $\phi_{(D_n,D_{n+1})}$ in the group $\text{Aut}(E^*\FF_0,U_{(D_n,D_{n+1})})$ such that 
$$\phi_{D_n}^{-1}\circ H_{D_n}^{-1}\circ H_{D_{n+1}}\quad\text{and}\quad\phi_{(D_n,D_{n+1})}$$
have same images in $\textup{Cent}(\text{Hol}_n)\left/_{<\text{Hol}_n>}\right.$. Hence, the automorphism
$$\phi_{D_n}^{-1}\circ H_{D_n}^{-1}\circ H_{D_{n+1}}\circ\phi_{(D_n,D_{n+1})}$$
acts in the local leaf. In the same way, the map $\phi_{(D_n,D_{n+1})}$ is an  automorphism of the restricted foliation in an open neighborhood of the set $U_{D_{n+1}}\cap U_{(D_n,D_{n+1})}$. Since the open set $U_{D_{n+1}}$ is conformally equivalent to a corona, there exists $\phi_{D_{n+1}}$ in $\text{Aut}(E^*\FF_0,U_{D_{n+1}})$ such that 
$$\phi_{D_{n+1}}^{-1}\circ \phi_{(D_n,D_{n+1})}$$ 
acts in the local leaf. Hence, we have obtained the induction hypothesis at rank $n+1$. 
\end{demo}
Let us now consider the following family of maps 
$$\left\{ \begin{array}{ll}
i \in \mathbb{I}_0= \textup{Comp}(\DD), & \psi_D=\phi_D\circ H_D\\
i \in \mathbb{I}_1=\textup{Comp}(\DD)^{2}, & \psi_{(D,D')}= \phi_{(D,D')}\circ H_D
\end{array}\right.$$
In view of the construction, the automorphism $\psi_i^{-1}\circ \psi_{j},\ i,j \in \mathbb{I}=\mathbb{I}_0\cup \mathbb{I}_1$ acts in the local leaf. In view of \cite{iso}, there exists a family of tangent vector fields $\fami{X_{ij}}_{i,j\in \mathbb{I}}$ and a family of holomorphic functions $\fami{t_{ij}}_{i,j\in\mathbb{I}}$ such that
$$\psi_i^{-1}\circ\psi_{j}=\Phi_{X_{ij}}^{t_{ij}},$$
where $\Phi_X^t$ refers to the flow of the vector field $X$ at $t$ time. Let us denote by $\mathcal{U}_i$ a neighborhood of $U_i$ in $\MM$. The deformation defined by the identification
$$s\longmapsto\left.\coprod_{i\in\mathbb{I}}{E^*\FF_0|}_{\mathcal{U}_i}\right/_{x\sim \Phi_{X_{ij}}^{s\cdot t_{ij}}x},\quad s\in \dbar$$
is well defined since the automorphisms of identification $x\rightarrow \Phi_{X_{ij}}^{s\cdot t_{ij}}x$ act in the local leaf. Moreover, in view of the cohomological interpretation of isoholonomic deformation, this deformation is precisely an isoholonomic deformation of foliation. The fiber of this deformation at $0\in \dbar$ is the foliation $\FF_0$ and the fiber at $1\in \dbar$ is a foliation analytically equivalent to $\FF_1$. Now, the theorem (\ref{exis}) applied with $C=\{0,1\}\subset \dbar$ ensures the existence of an isoholonomic deformation $\mathcal{R}$ over $\dbar$ from $\FF_0$ to $\FF_1$ such that the underlying deformation of separatrix is the trivial deformation $S\times \dbar$. 
Since $\FF_0$ is d-quasihomogeneous, the deformation $\mathcal{R}$ is locally trivial. Hence, $\FF_0$ and $\FF_1$ are analytically conjugated. This ends the proof in the convergent context. The transposition to the formal context is let to the reader.

{\small 
\bibliography{biblithese}
\bibliographystyle{plain}}
\end{document}